\documentclass[10pt]{amsart}
\usepackage{amsmath,amsthm,amssymb,latexsym,graphicx,verbatim,enumerate,%
ifpdf,mdwlist}
\ifpdf
\usepackage{hyperref}
\else
\usepackage[hypertex]{hyperref}
\fi

\title[$\Sigma^1_1$-Completeness of Decomposable Groups]{The Decomposability
Problem for Torsion-Free Abelian Groups is Analytic-Complete}

\author{Kyle Riggs}
\address{Department of Mathematics\\
Indiana University\\
Bloomington, IN 47405\\
USA}
\email{\href{mailto:kwriggs@umail.iu.edu}{kwriggs@umail.iu.edu}}
\urladdr{\url{http://www.math.indiana.edu/~kwriggs/}}


\thanks{I would like to thank Steffen Lempp and Alexander Melnikov for their
fruitful discussions with me.}

\begin{document}

\begin{abstract}
We discuss the decomposability of torsion-free abelian groups. We show that
among computable groups of finite rank this property is $\Sigma^0_3$-complete. However,
when we consider computable groups of infinite rank, it becomes $\Sigma_1^1$-complete
(and $\boldsymbol{\Sigma^1_1}$-complete for groups of infinite rank in general),
so it cannot be characterized by a first-order formula in the language of arithmetic.
\end{abstract}

\maketitle

\section{Introduction}

A group is computable if its domain can be enumerated effectively and the binary
operation of the group is computable (via the enumeration). Such an enumeration is 
called a \textit{computable presentation}. In other words, a computable group
has a computable word problem. As with many other computable structures,
there are a variety of questions we can ask about them. Perhaps
the first question that comes to mind is what groups can be presented computably.
This type of question is typically answered by looking at a specific class
of groups and identifying which groups in that class have computable copies.
Downey and others~\cite{DC} studied the class of completely decomposable
groups.
\newline

\textbf{Definition 1.1:} An abelian group $G$ is \textit{completely decomposable}
if it can be written
$$G = \underset{i}{\oplus} H_i$$
where each $H_i$ is a subgroup of $Q$.
\newline

In particular, they considered groups of the form
$$G_S = \underset{p \in S}{\oplus} Q_p,$$
where $S$ is a set of primes and $Q_p$ is the subgroup of $Q$ generated by
the elements $\{\frac{1}{p^k}: k \in \omega\}$. They found that there is a computable
copy of $G_S$ iff $S$ is $\Sigma^0_3$.

Khisamiev~\cite{Khis} looked at (countable) reduced torsion groups,
which are uniquely determined by their Ulm sequences. He was able to
completely characterize the Ulm sequences of length $< \omega^2$ which can
occur in a computably presented group. Ash, Knight, and Oates~\cite{AKO}, working slightly
later, independently duplicated his results.
\newline

Another type of question that arises involves the notion of computable
categoricity.

\textbf{Definition 1.2:} A computably presentable group is \textit{computably
categorical} if any two computable copies of the group have a computable
isomorphism between them.

Additionally, a computably presentable group is $\Delta^0_n$-\textit{categorical}
if any two computable copies of the group have a $\Delta^0_n$-isomorphism
between them.
\newline

Downey and Melnikov~(\cite{EC},\cite{DM}) studied homogeneous completely
decomposable groups, in which each summand has the same type (see Definition
2.2). They found that every completely decomposable group is $\Delta^0_5$
-categorical, but homogenous completely decomposable groups are actually
$\Delta^0_3$-categorical.
\newline

Studying computable groups can also yield results telling us the difficulty
of determining some of their most fundamental properties. For example, Downey
and Montalban~\cite{IP} studied the isomorphism problem for torsion-free abelian
groups. They found that:
\newline

1) the set of pairs of computable indices for isomorphic torsion-free abelian
groups is $\Sigma^1_1$-complete, and

2) the set of isomorphic pairs of torsion-free abelian groups is $\boldsymbol{\Sigma^1_1}$
-complete
\newline

In their research of completely decomposable groups, Downey and Melnikov~\cite{DM}
found that the index set of completely decomposable groups can be described
by a $\Sigma^0_7$ formula (though it is not known if this is sharp). The
property with which this paper is concerned is a similar one:~decomposability.
\newline

\textbf{Definition 1.3:} An abelian group is \textit{decomposable} if it can
be written as the direct sum of two (or more) nontrivial subgroups. Otherwise,
it is \textit{indecomposable}.
\newline

Often the best way to study an abelian group is by writing it as a direct
sum of its indecomposable subgroups, so determining whether a group is
decomposable is a problem at the heart of abelian group theory.

It is known that the only indecomposable torsion groups are the cocyclic groups (groups
of the form $Z(p^k)$, with $k \in \omega + 1$), and that every mixed group
is decomposable. Torsion-free groups of rank 1 are indecomposable, but
beyond this no classification has been found. It has been conjectured by some
(Kudinov, Melnikov) that this is in part because the class of torsion-free
decomposable groups is not arithmetical. We shall see that this is indeed
the case. In fact, the class is not even hyperarithmetical.

After a short discussion of some concepts from logic and algebra, we will construct an
indecomposable group of rank 2. We will apply the ideas used in this example
to show that the index set of computable decomposable groups of finite rank is
$\Sigma^0_3$-complete. We will then shift our focus to groups of infinite
rank and show that:
\newline

1) the index set of computable decomposable torsion-free abelian groups is
$\Sigma^1_1$-complete, and

2) the set of decomposable torsion-free abelian groups is $\boldsymbol{\Sigma^1_1}$
-complete
\newline

\section{Some Computational Complexity Hierarchies}

The arithmetical hierarchy was developed to describe the complexity of properties based on their formulas in the language of arithmetic. A computable set (or relation) $S \subset \omega$ is said to be $\Sigma^0_0$ (or $\Pi^0_0$). A set $S_1$ is $\Sigma^0_{n+1}$ if it can be characterized by a formula of the form
$$x \in S_1 \Leftrightarrow (\exists y \in \omega) \, R_1(x,y)$$
where $R_1$ is a $\Pi^0_n$ relation. Likewise, a set $S_2$ is $\Pi^0_{n+1}$ if it can be characterized by a formula of the form 
$$x \in S_2 \Leftrightarrow (\forall y \in \omega) \, R_2 (x,y)$$
where $R_2$ is a $\Sigma^0_n$ relation.

In other words, the $n$ represents how many times the formula alternates quantifiers over $\omega$ (or some other infinite computable set), and a $\Sigma^0_n$ formula starts with an existential quantifier, while a $\Pi^0_n$ formula starts with a universal quantifier.

For example, given a computable group $G$, an element $g \in G$, and a fixed prime $p$, there is a $\Pi^0_2$ formula that says whether $p$ infinitely divides $g$
$$p \, |^\infty \, g \Leftrightarrow (\forall k \in \omega) \; (\exists h \in G) \; p^k h = g$$
\newline

Any set that is characterized by a $\Sigma^0_n$ or $\Pi^0_n$ formula for some $n$ is said to be \textit{arithmetical}. If we allow quantifiers over functions from $\omega$ to $\omega$ (or between any two computable sets), then our formula will be \textit{analytic}. We say that a formula is $\Sigma^1_1$ if it is of the form
$$(\exists f \in \omega^\omega) \; R(f)$$
where $R$ is any arithmetical formula. Similarly, a $\Pi^1_1$ formula is of the form
$$(\forall f \in \omega^\omega) \; R(f)$$
where $R$ is arithmetical.

For any complexity class $\Gamma$, we say that a set $A$ is $\Gamma$-complete if any other set $B \in \Gamma$ can be "coded" into $A$. That is to say, there is a computable function $f: \omega \rightarrow \omega$ such that $x \in B$ iff $f(x) \in A$.

The canonical example of a $\Sigma^1_1$-complete set is the index set of computable trees in $\omega^{<\omega}$ with an infinite path.

The Borel hierarchy is a complexity class structure for Polish spaces (like $\omega^\omega$) which defines $\boldsymbol{\Sigma^0_1}$ sets to be open sets and $\boldsymbol{\Pi^0_1}$ sets to be closed sets. In this hierarchy, a set is $\boldsymbol{\Sigma^1_1}$ or \textit{analytic} if it is the image of a Polish space under a continuous mapping. For example, the set of trees (not necessarily computable) in $\omega^{<\omega}$ with an infinite path is $\boldsymbol{\Sigma^1_1}$-complete.
\newline

To summarize, this paper shows that 

1) the index set of computable decomposable torsion-free abelian groups of finite rank is $\Sigma^0_3$-complete,

2) the index set of computable decomposable torsion-free abelian groups of infinite rank is $\Sigma^1_1$-complete, and

3) the set of decomposable torsion-free abelian groups of infinite rank is $\boldsymbol{\Sigma^1_1}$-complete.

\section{Algebra Background}

In this paper, we will exclusively discuss torsion-free abelian groups. Also,
the term \textit{basis} will refer to a maximal linearly independent subset
of a group. The rank of a group is the cardinality of any of its bases.
\newline

\textbf{Definition 3.1:} Given an abelian group G, an element $x \in G$, and a prime $p$, the \textit{height} of $p$ at
$x$ is given by
$$h_p(x) = \sup\{k : p^k | x\}$$
We call
$$\chi_G(x) = (h_2(x),h_3(x), h_5(x),...)$$
the \textit{characteristic} of $x$ in $G$.
\newline

\textbf{Definition 3.2:} We define an equivalence relation on characteristics by saying that 
$\chi_G(x)~\sim~\chi_G(y)$ if
\begin{itemize}
\item for all $p, \; h_p(x) = \infty \Leftrightarrow h_p(y) = \infty,$ and
\item $h_p(x) = h_p(y)$ for all but finitely many $p$
\end{itemize}

We call the equivalence classes \textit{types}. In other words, $x$ and $y$
have the same type iff there exist integers $m$ and $n$ such that
$\chi_G(mx) = \chi_G(ny)$.
\newline

We can put a partial order on types by declaring for two types $\alpha, \beta$
that $\alpha \preceq \beta$ if, given any element $a$ of type $\alpha$ and any
element $b$ of type $\beta$,
\begin{itemize}
\item for all $p, \; h_p(a) = \infty \Rightarrow h_p(b) = \infty$ and
\item $h_p(a) \leq h_p(b)$ for all but finitely many $p$
\end{itemize}
\vspace{.1 in}

\textbf{Definition 3.3:} A nonzero element has \textit{strictly maximal type}
if no nonzero element linearly independent from it has a greater or equal type.
\newline

There is a class of subgroups which we must introduce before giving any proofs.
\newline

\textbf{Definition 3.4:} In an abelian group $G$, a subgroup $H$ is called
\textit{pure} if for every $x \in H$ and $m \in \omega$, if $m$ divides $x$ in $G$,
$m$ also divides $x$ in $H$.

If $S$ is a set of elements in $G$, the \textit{pure subgroup generated by}
 $S$ is the smallest pure subgroup of $G$ containing $S$.
\newline

For example, $Z$ is a pure subgroup of $Z^2$, but $Z$ is not a pure subgroup
of $Q$. In fact, the only pure subgroups of $Q$ are $Q$ and the trivial
subgroup.
\newline

\section{An Example of an Indecomposable Group} The following example can be
found in Fuchs~\cite{Fuchs}. Let $G_0$ be the free abelian group generated by two
elements, $x_1$ and $x_2$. For every $k > 0$, we add elements of the form
$$\frac{x_1}{3^k} \text{ and } \frac{x_2}{5^k}$$
\newline
to $G_0$. We also add the element $\frac{x_1 + x_2}{2}$ to the group. 
We denote by $G$ the group generated by all these elements.

Note that $\{x_1,x_2\}$ is still a basis for this group, and that 
$$\chi_G(x_1) = (0,\infty,0,0,0,...) \text{ and } \chi_G(x_2) = (0,0,\infty,0,0,...)$$
Furthermore, any element of the form $q_1 x_1 + q_2 x_2$ with both coefficients
nonzero has type $(0,0,0,...)$. Thus, $x_1$ and $x_2$ both have strictly
maximal type.
\newline

\textbf{Proposition 4.1:} In a decomposable group $G \; ( = A \oplus B)$, if
$x \in G$ decomposes as $x = a + b$ and an integer $m$ divides $x$, then $m$
divides $a$ and $b$ as well.
\newline

\textit{Proof}: Let $y \in G$ be such that $my=x$, and suppose $y$ decomposes
$y = a_1 + b_1$. Then we see that
$$ma_1 + mb_1 = my =  x = a + b$$
$ma_1 \in A$ and $mb_1 \in B$, so $ma_1 = a$ and $mb_1 = b$.
\qed \newline

\textbf{Corollary 4.2:} In a decomposable group $G$, every element of strictly
maximal type must be contained in a direct summand.
\newline

\textit{Proof}: Suppose that $x$ is an element of strictly maximal type that
is not in $A$ or $B$. Then we can write $x = a + b$ with $a \in A$ and $b \in
B$, and both $a$ and $b$ nonzero. Because $x$ has strictly maximal type, there is a
prime $p$ and an integer $k$ such that $p^k$ divides $x$, but neither $a$ nor
$b$. However, this contradicts the proposition.
\qed \newline

Thus, we can assume $x_1 \in A$ and $x_2 \in B$. Now consider the decomposition
of the element
$$\frac{x_1 + x_2}{2} = a + b$$
\newline
with $a \in A$ and $b \in B$. It is clear that $2a = x_1$ and $2b = x_2$, but
there are no elements in $G$ which satisfy these equations. Thus, the group
is indecomposable.
\newline

The proofs contained in this paper will mimic this technique of creating
elements of strictly maximal type and then introducing elements which force
them to be contained in the same direct summand. We call these elements
\textit{links}.
\newline

\textbf{Definition:} Let $x$ and $y$ be two elements of strictly maximal type
in a torsion-free abelian group $G$. If there is a prime $p$ which divides
the sum $x+y$ but neither $x$ nor $y$, then the element $\frac{x+y}{p}$ is
a link connecting $x$ and $y$. We say that $x$ and $y$ are connected by a
\textit{chain} of links if there are elements $x_1, x_2, ..., x_n$ such that
the sequence $\{x_0 = x, x_1, x_2, ..., x_n, x_{n+1} = y\}$ has the property
that for $0 \leq i \leq n$, there is a link connecting $x_i$ and $x_{i+1}$.
\newline

The following proposition gives us a simple way to construct indecomposable
groups.
\newline

\textbf{Proposition 4.3:} If a torsion-free group has a basis of elements of
strictly maximal type, with each pair of them having a link or a chain of
links connecting them, then it is indecomposable.
\newline

\textit{Proof}: Every element of strictly maximal type must be contained in
a direct summand, and any two of these elements with a link between them
must be in the same direct summand. Transitively, this is also true of any
two elements with a chain of links connecting them. Thus, the entire basis
is contained in a single direct summand, so the group is indecomposable.

\section{Groups of Finite Rank}

\textbf{Remark 5.1:} Let $G$ be a group of finite rank, and assume $G = A
\oplus B$. Then
\begin{enumerate} 
\item $\text{rank}(G) = \text{rank}(A) + \text{rank}(B)$
\item If $\{a_1,...,a_n\}$ is a basis for $A$ and $\{b_1,...,b_m\}$ is a 
basis for $B$, then \\ $\{a_1,...,a_n,b_1,...,b_m\}$ is a basis for
$G$ with the following property:
\end{enumerate}

\begin{align*}
\text{If there exists an element } g = \sum\limits_{i=1}^n q_i a_i + 
\sum\limits_{j=1}^m r_j b_j, \text{ then there exist} \\
\text{elements } g_A, g_B \text{ such that } g_A =
\sum\limits_{i=1}^n q_i a_i \text { and } g_B = \sum\limits_{j=1}^m r_j b_j
\end{align*}

Conversely, if $\{a_1,...,a_n,b_1,...b_m\}$ is a basis for $G$ with this
property, then the pure subgroup generated by the $a_i$'s and the pure 
subgroup generated by the $b_j$'s give a decomposition of $G$. Thus, a group
of finite rank is decomposable iff it has a basis with this property.

Let $[G]^{<\omega}$ denote the set of all finite sets in $G$. The following
$\Pi^0_2$ formula describes the property of being a basis of $G$. For
$\bar{x} \in [G]^{<\omega}$,
\begin{align*}
 BASIS(\bar{x}) \Leftrightarrow [(\forall y \in G) \; (\exists \bar{q} 
\in Q^{<\omega}) \; (|\bar{q}| = |\bar{x}| \wedge y = \sum\limits_i q_i x_i) \\
\wedge (\forall  \bar{q} \in Q^{<\omega})\; (|\bar{q}| = |\bar{x}| \wedge
\sum\limits_i q_i x_i = 0) \Rightarrow \bar{q} = \bar{0}]
\end{align*}

If we take the conjunction of that formula with one describing the property
in Remark 4.1, we have the following $\Sigma^0_3$ formula for decomposable
groups of finite rank:

\begin{align*}
(\exists \, \bar{a},\bar{b} \in [G]^{<\omega}) \; \{BASIS(\bar{a} \sqcup
\bar{b}) \wedge \bar{a} \neq \o \wedge \bar{b} \neq \o \wedge (\forall y \in
G)\, (\forall \bar{q} \in Q^{<\omega})\\
(\exists w \in G) [(|\bar{q}| = |\bar{a}| + |\bar{b}| \wedge y = \sum\limits_i
q_i a_i + \sum\limits_j q_j b_j) \Rightarrow \; w = \sum\limits_i q_i a_i]\}
\end{align*}
\newline

\textbf{Theorem 5.2:} The set of decomposable groups of finite rank is
$\Sigma^0_3$-complete.
\newline

\textit{Proof}: Recall that $Cof = \{n : W_n \text{ is cofinite}\}$ is
$\Sigma^0_3$-complete. In order to prove our result, we construct a function
from $\omega$ to groups of rank 2 such that $G_n$ is decomposable iff $W_n$ is
cofinite.
\newline

\textit{Construction}: We start with a group $G$ generated by the following elements:

$$ \langle g_1, g_2, \frac{g_1 + g_2}{2}, \frac{g_1}{3}, \frac{g_2}{5}, 
\frac{g_1}{7},\frac{g_2}{11},... \rangle $$
\newline
($g_1$ and $g_2$ are linearly independent). 

$g_1$ is divisible by all odd-indexed primes, and $g_2$ is divisible by all
even-indexed primes (except $p_0=2$), so they have incomparable (indeed,
strictly maximal) types. Thus, like the example above, our initial group
$G$ is indecomposable.

The group $G_n$ is generated by adding $\frac{g_2}{p_{2k+1}}$ for every $k$
such that $\Phi_n(k) \downarrow$.
\newline

\textit{Verification}: If $W_n$ is coinfinite, then $g_1$ is still divisible
by infinitely many primes that do not divide $g_2$. Thus, the types remain
incomparable, and the group remains indecomposable.
\newline

If $W_n$ is cofinite, then the type of $g_2$ is strictly greater than the type
of $g_1$. There are finitely many primes that divide $g_1$ but not $g_2$.
Denote their product by $m$.
\newline

\textbf{Lemma 5.3:} $G_n = A \oplus B$, where $A$ is the pure subgroup generated
by $a = \frac{g_1 + m g_2}{2}$ and $B$ the pure subgroup generated by $g_2$.
\newline

\textit{Proof}: We observe that $\frac{g_1 + g_2}{2} = a - \frac{m-1}{2} g_2$ (Note that
$m$ is a product of odd primes).

Any element of the form $\frac{g_1}{p} \in G_n$ can be written
$$\frac{g_1}{p} = \frac{2}{p}a - \frac{m}{p}g_2$$

If $p \nmid m$, then $p \, | \, g_2$, so $\frac{m}{p}g_2 \in B$. Thus, every
generating element of the group can be uniquely decomposed, so the group is
decomposable.
\qed\newline

$G_n$ is decomposable iff $W_n$ is cofinite, so the theorem is proved. \qed\newline

\section{Groups of Infinite Rank}

We can adapt the formula used for groups of finite rank to describe decomposable
groups of infinite rank. However, this means the first existential quantifier
is searching over infinite sets instead of finite sets, so the $\Sigma^0_3$
formula becomes a $\Sigma^1_1$ formula (here $BASIS$ is a $\Pi^0_2$-formula
on infinite sets):

\begin{align*}
(\exists \, \bar{a},\bar{b} \in [G]^{\leq\omega}) \; [BASIS(\bar{a} \sqcup
\bar{b}) \wedge \bar{a} \neq \o \wedge \bar{b} \neq \o \wedge (\forall y 
\in G)\\ (\forall \bar{q} \in Q^{<\omega}) (\exists w \in G) (y =
\sum\limits_i q_i a_i + \sum\limits_j q_j b_j) \Rightarrow \; w = \sum\limits_i q_i a_i)]
\end{align*}
\newline

\textbf{Theorem 6.1:} The set of decomposable groups of infinite rank is
$\boldsymbol{\Sigma^1_1}$-complete.
\newline

\textit{Proof}: We will construct a function from trees in $\omega^{<\omega}$
to torsion-free abelian groups of infinite rank that takes a tree $T$ and
gives a group $G_T$ that is decomposable iff $T$ has an infinite path. (Recall
that the set of trees in $\omega^{<\omega}$ which have an infinite path is
$\boldsymbol{\Sigma^1_1}$-complete.)
\newline

\textit{Construction of the group $G$}: We start with a countably infinite set
of linearly independent elements: $x_1, x_2,...$ and $\{x_\sigma\}_{\sigma
\in \omega^{< \omega}}$ (which we denote as the $x$-elements), and $y_1,
y_2,...$ (the $y$-elements). These elements form a basis for our group. We
will give them all strictly maximal type and introduce links connecting all
the $x$-elements and separate links connecting all the $y$-elements.
\newline

The initial group $G_0$ is generated by the following elements:

\begin{itemize}

\item For $i,k > 0$ and $\sigma \in \omega^{<\omega}$, 
$$\frac{x_i}{p_{\langle 0,i \rangle}^k} \, , \; \frac{y_i}{p_{\langle 1,i \rangle}^k}
\text{ , and } \frac{x_\sigma}{p_{\langle 2,\sigma \rangle}^k}$$

\item For $0 < i<j$, $$\frac{x_i + x_j}{p_{\langle 3,\langle i,j \rangle \rangle}}
\text{ and } \frac{y_i + y_j}{p_{\langle 4, \langle i,j \rangle \rangle}}$$

\item For every $i \geq 0$ and $\sigma, \rho \in \omega^{<\omega}$, 
$$\frac{x_i + x_\sigma}{p_{\langle 5, \langle i,\sigma\rangle \rangle}}
\text{ and } \frac{x_\sigma + x_\rho}{p_{\langle 6, \langle \sigma,\rho
\rangle \rangle}}$$

\item For $n > 1$, $$\frac{y_1 + y_2 + ... + y_n}{p_{\langle 7,n \rangle}}$$
\end{itemize}

All the $x$- and $y$-elements are elements of strictly maximal type, and due
to the links, all the $x$-elements must be in the same direct summand of 
$G_0$ (as do the $y$-elements). Thus, $G_0$ can only be decomposed as $G_0
= A \oplus B$, where $A$ is the pure subgroup containing all the $x$-elements,
and $B$ is the pure subgroup containing all the $y$-elements.
\newline

Now we add to $G_0$ links of the form $$\frac{x_i + y_i}{p_{\langle 8,
i \rangle}}$$ for $i \geq 0$, and denote by $G$ the group generated by these
elements. Now every $x$- and $y$-element are connected by a chain of links,
so $G$ is indecomposable.
\newline

\textit{Construction of $G_T$}: Given a tree $T$ in $\omega^{<\omega}$, we
will add elements to $G$ to form a group $G_T$ that will be decomposable
iff $T$ has an infinite path through it. The idea is that if there is an
infinite path $\pi$, then $G_T = A_T \oplus B_\pi$, where $A_T$ is the pure
subgroup containing the $x$-elements, and $B_\pi$ is the pure
subgroup of $G_T$ containing all the elements of the form $y_i + x_{\pi
\upharpoonright i}$. Note that if there is more than one infinite path through
$T$, there will be more than one way to decompose $G_T$.
\newline

Enumerate $T$ so that each string in $T$ is enumerated after all of its initial
segments. When we see $\sigma \in T$ with $|\sigma| = n$, we do the following:
(It's worth noting that in each case, the introduction of the first element
creates the second element. We list both simply to remind the reader that
the second element also exists)
\newline

(1) For $i \leq n$, we add to the group the elements  $$\frac{y_i + x_{\sigma
\upharpoonright i}}{p_{\langle 1,i \rangle}^n} \text{ and } \frac{x_{\sigma
\upharpoonright i}}{p_{\langle 1,i \rangle}^n}$$

(2) For $i < n$, we add to the group the elements $$\frac{(y_i + x_{\sigma
\upharpoonright i}) + (y_n + x_\sigma)}{p_{\langle 4, \langle i,n \rangle
\rangle}} \text{ and } \frac{x_{\sigma \upharpoonright i} + x_\sigma}{p_{\langle
4, \langle i,n \rangle \rangle}}$$

(3) We add to the group the elements $$\frac{y_n + x_\sigma}{p_{\langle 8,n
\rangle}} \text{ and } \frac{x_n - x_\sigma}{p_{\langle 8,n \rangle}}$$

(4) Finally, we add the elements $$\frac{(y_1 + x_{\sigma \upharpoonright 1})
+ (y_2 + x_{\sigma \upharpoonright 2}) + ... + (y_n + x_\sigma)}{p_{\langle
7,n \rangle}} \text{ and } \frac{x_{\sigma \upharpoonright 1} + x_{\sigma
\upharpoonright 2} + ... + x_\sigma}{p_{\langle 7,n \rangle}}$$
\newline

\textit{Verification}: If an infinite path $\pi$ does exist, we shall see that
$G_T = A_T \oplus B_\pi$ (as described above). 

Each $x_i$ and $x_\sigma$ is contained in $A_T$.  We have $y_j = -x_{\pi
\upharpoonright j} + (y_j + x_{\pi \upharpoonright j})$. Both of these elements
are infinitely divisible by $p_{\langle 1,j \rangle}$ because $x_{\pi
\upharpoonright j}$ went through step (1) infinitely often.
\newline

For $0 < i < j$,  $$\frac{y_i + y_j}{p_{\langle 4, \langle i,j \rangle \rangle}}
= \frac{(y_i + x_{\pi \upharpoonright i}) + (y_j + x_{\pi \upharpoonright j})}
{p_{\langle 4, \langle i,j \rangle \rangle}} - \frac{x_{\pi \upharpoonright i}
+ x_{\pi \upharpoonright j}}{p_{\langle 4, \langle i,j \rangle \rangle}}$$
These elements were created during step (2) of some stage.
\newline

For $i > 0$, $$\frac{x_i + y_i}{p_{\langle 8,i \rangle}} = \frac{x_i - x_{\pi
\upharpoonright i}}{p_{\langle 8,i \rangle}} + \frac{y_i + x_{\pi \upharpoonright
i}}{p_{\langle 8,i \rangle}}$$
These elements were created during step (3) of some stage.
\newline

For $n > 1$, 
\begin{multline*}
\frac{y_1 + y_2 + ... + y_n}{p_{\langle 7,n \rangle}} =
-\frac{x_{\pi \upharpoonright 1} + x_{\pi \upharpoonright 2} + ... + x_{\pi
\upharpoonright n}}{p_{\langle 7,n \rangle}} +{} \\
\frac{(y_1 + x_{\pi \upharpoonright
1}) + (y_2 + x_{\pi \upharpoonright 2}) +  ... + (y_n + x_{\pi \upharpoonright
n})}{p_{\langle 7,n \rangle}}
\end{multline*}

These elements were created during step (4) of some stage.

We see that all the generating elements of $G_T$ can be uniquely decomposed,
so $G_T = A_T \oplus B_\pi$.
\newline

Now suppose $G_T$ is decomposable as $G_T = A' \oplus B'$. All the$x$-elements
still have strictly maximal type, so they must be in the
same direct summand ($A'$). 

Each $y_j$ can be decomposed $y_j = a_j + b_j$, where $a_j \in A'$ and
$b_j \in B'$. We know $a_j$ and $b_j$ are infinitely divisible by $p_{\langle
1,j \rangle}$ because $y_j$ is. The only other basis elements that could be
infinitely divisible by this prime are the elements $x_\sigma$ with $|\sigma| = j$.
\newline

\textbf{Lemma 6.2:} If there is any $y_j \in A'$, then $G_T = A'$ (and
$B' = 0$).
\newline

\textit{Proof}: Suppose $y_j \in A' \; (y_j = a_j)$, and that another element
$y_i \notin A'$. Let $q = p_{\langle 4, \langle j,i \rangle \rangle}$ (we can assume
that $j < i$). Then $q$ must divide $y_i + y_j$,
and thus, $b_i + b_j$ (which is just $b_i$). There is some finite sum such that

$$k_0 b_i = k_1 y_i + \sum\limits_{|\sigma|=i} k_\sigma x_\sigma$$
with each $k_\sigma,k_0,k_1 \in \omega$. (Recall that the only other basis elements
that could be infinitely divisible by $p_{\langle 1, i \rangle}$ are the
elements $x_\sigma$ with $|\sigma| = i$.)  We can also write

$$k_0 a_i = (k_0 - k_1) y_i - \sum\limits_{|\sigma|=i} k_\sigma x_\sigma$$
Note that if $k_0 \neq k_1$, then $y_i \in A'$. Thus, $k_1 = k_0$, so

$$b_i = y_i + \frac{1}{k_0} \sum\limits_{|\sigma| = i} k_\sigma x_\sigma$$
and this must be divisible by $q$. However, $q$ does not divide $y_i$, nor
any nontrivial linear combination of $y_i$ with $x$-elements (though $q$ 
does divide $y_i + y_j$). Therefore, $y_i \in A'$. This is true for every
$y_i$, so $A' = G_T$.
\qed\newline

There is no $y_j \in B'$, either. This is because there are no elements
$$\frac{x_j} {p_{<8,j>}} \, , \, \frac{y_j}{p_{<8,j>}}$$
So we see that every $y_j = a_j + b_j$, with both components being nonzero.
\newline

Now suppose $y_1, y_2$ decompose as $$y_1 = \sum\limits_{|\sigma| = 1} k_\sigma
x_\sigma + b_1 \text{ and } y_2 = \sum\limits_{|\rho| = 2} l_\rho x_\rho + b_2$$

We shall denote $p_{\langle 7,2 \rangle}$ by $r$. $r | (y_1 + y_2)$, so it
must also divide
$$a_1 + a_2 = \sum\limits_{|\sigma| = 1} k_\sigma x_\sigma
+ \sum\limits_{|\rho| = 2} l_\rho x_\rho$$
Although $r$ does not divide any $x_\sigma$, from step (4) of the construction
we see that $r$ divides elements of the form $x_\sigma + x_{\sigma\, \hat{ }
\, m}$ where $|\sigma| = 1$ (and $\sigma \, \hat \, m \in T$). Thus, $r$ also
divides elements of the form
$$\sum\limits_m l_m (x_\sigma + x_ {\sigma \, \hat{} \, m}) = kx_\sigma +
\sum\limits_m l_m x_{\sigma\, \hat{ } \, m}$$
where $k = \sum\limits_m l_m$ (and $\sigma \, \hat{ } \, m \in T$).
\newline

From this we see that
$$r | \; (\sum\limits_{|\sigma| = 1} k_\sigma x_\sigma + \sum\limits_{|\rho|
= 2} l_\rho x_\rho) \text{ iff } \, k_\sigma \equiv \sum\limits_{\rho \succ
\sigma} l_\rho (\text{mod } r)$$
for each $\sigma$ with $|\sigma| = 1$.

Similarly, $p_{\langle 7,3 \rangle} \, | \, (y_1 + y_2 + y_3)$, so for each
$\tau \in T$ with $|\tau| = 3, \; p_{\langle 7,3 \rangle} \, | \, (x_\sigma
+ x_\rho + x_\tau)$, where $\sigma \prec \rho \prec \tau$.
\newline

By the same reasoning, we see that if $y_3 = \sum\limits_{|\tau| = 3} m_\tau
x_\tau + b_3$, then for each $\sigma \in T$ with $|\sigma| = 1$,

$$k_\sigma \equiv \sum\limits_{\rho \succ \sigma} l_\rho \equiv 
\sum\limits_{\tau \succ \sigma} m_\tau \; (\text{mod } p_{\langle 7,3 \rangle})$$

There are infinitely many such equivalences, so we see that 

$$k_\sigma = \sum\limits_{\rho \succ \sigma} l_\rho = \sum\limits_{\tau \succ
\sigma} m_\tau = ...$$

It is also true that for each $\rho \in T$ with $|\rho| = 2$,

$$l_\rho \equiv \sum\limits_{\tau \succ \rho} m_\tau \;
(\text{mod } p_{\langle 7,3 \rangle})$$

Continuing this process, we see that the following also holds:

$$l_\rho = \sum\limits_{\tau \succ \rho} m_\tau = ...$$

Thus, if we choose a $\sigma$ of length 1 such that $k_\sigma \neq 0$ (which
we are guaranteed by the fact that $y_1 \notin B'$), there must be a $\rho$
of length 2 such that $\sigma \prec \rho$ and $\l_\rho \neq 0$, and a $\tau$
of length 3 such that $\rho \prec \tau$ and $m_\tau \neq 0$. By repeating
this process, we find an infinite path through $T$.

Thus, $G_T$ is decomposable iff $T$ has an infinite path.
\qed
\newline


\bibliographystyle{plain}

\end{document}